\documentclass{amsart}
\usepackage{amssymb,amsmath,xcolor,graphicx,xspace,colortbl, rotating}
\usepackage{amsfonts}
\usepackage{amsmath}
\usepackage{amssymb}
\usepackage[T1]{fontenc}
\usepackage{graphics}
\usepackage{graphicx}
\usepackage{bm}

\DeclareGraphicsExtensions{.pdf,.svg,.eps,.ps,.png,.jpg,.jpeg}
\newtheorem{theorem}{Theorem}
\newtheorem{lemma}[theorem]{Lemma}

\newtheorem{definition}[theorem]{Definition}

\begin{document}
\title{Counting words satisfying the rhythmic oddity property}
\author{Franck Jedrzejewski}
\maketitle

\begin{abstract}
This paper describes an enumeration of all words having a combinatoric
property called ``rhythmic oddity property''named \emph{rop-words}.\ This
property was introduced by Simha Aron in the 1990s. The set of rop-words is
not a subset of the set of Lyndon words, but is very closed. We show that
there is a bijection between some necklaces and rop-words. This leads to a
formula for counting the rop-words of a given length. \textsc{Keywords:}
Combinatoric on words. Lyndon words. Rhythmic oddity. Music formalization
\end{abstract}


\bigskip

The \emph{rhythmic oddity property}\ was discovered by ethnomusicologist
Simha Aron \cite{Aro1991} in the study of Aka pygmies music.\ The rhythms
satisfying this property are refinement of aksak rhythms described by C.\ Br%
\u{a}iloiu \cite{Bra1952} in 1952, and used also in Turkish and Bulgarian
music. They have been studied by M.\ Chemillier and C.\ Truchet in \cite%
{Che2003} and Andr\'{e} Bouchet in \cite{Bou2010} who gave some nice
characterizations. Related problems of asymmetric rhythms have been studied
by Rachael Hall and Paul Klingsberg \cite{Hal2004, Hal2006}. In this paper,
we carry on these studies by showing a one to one correspondance between
some necklaces and words satisfying the rhythmic oddity property.

\section{The rhythmic oddity property}

Patterns with rhythmic oddity property are combinations of durations equal
to 2 or 3\ units, such as the famous Aka pygmies rhythm 32222322222, and
such that when placing the sequence on a circle, ``one cannot break the
circle into two parts of equal length whatever the chosen breaking point.''\
In the langage of combinatorics of words, this property in terms of words $%
\omega $ over the alphabet $A =\{2 ,3\}$ is defined as follows.\ The height $%
h (\omega )$ of a word $\omega =\omega _{0} \omega _{1} \ldots \omega _{n
-1} $ of length $n$ is by definition the sum of its letters $h (\omega ) =
\sum _{j =0 ,\ldots ,n -1}\omega _{j}$.\ A word $\omega $ satisfies the
rhythmic oddity property (rop) if $h (\omega )$ is even and no cyclic shift
of $\omega $ can be factorized into two words $u v$ such that $h (u) =h (v)%
\text{.}$ For short, we call \emph{rop}-word a word over the alphabet $\{2
,3\}$ satisfying the rhythmic oddity property.\ For instance, the word $%
32322 $ of height 12 is a rop-word, as well as all words of the form $32^{n}
32^{n +1}$ for all non-negative integers $n$, where the notation $2^{n}$
means the letter 2 is repeated $n$ times.\ The properties of rop-words have
been outine in \cite{Che2003}.\ But contrary to what is sometimes written,
the set of rop-words is not a subset of the set of Lyndon words.\ The words
222 and 233233233 are rop-words, but not Lyndon words. A Lyndon word is a
string that is strictly smaller in lexicographic order than all of its
rotations.\ Conversely, the set of Lyndon words is not included in the set
of rop-words, since 2233 is a Lyndon word, but not a rop-word (the words 23
and 32 have the same height and 2332 is a rotation of 2233). We call a
Lyndon rop-word a word of the monoid $\{2 ,3\}^{ \ast }$ that is both a
Lyndon word and a rop-word. The aim of this paper is to count the number of
Lyndon rop-words and the number of rop-words of length $n$.

\section{A\ bijection between some necklaces and rop-words}

Andr\'{e} Bouchet gave some characterizations of rop-words in \cite{Bou2010}%
.\ We present the main results of his paper.\ Let $\varepsilon $ be the
empty word and $\omega $ a word over $\{2 ,3\}$. The cycle $\delta $\ of $%
\omega $ is defined by $\delta (\varepsilon ) =\varepsilon $ and $\delta (a
\omega ) =\omega a\text{,}$ for $a \in \{2 ,3\}$. The rotations of $\omega $%
\ are the words $\delta ^{k} (\omega )\text{,}$ for all positive integer $k
>0$. In his paper, Bouchet shows the following lemma.

\begin{lemma}
Let $\omega =\omega _{0} \omega _{1} \ldots \omega _{n -1}$ be a word over
the alphabet\ $A =\{2 ,3\}$ of height $2 h$.\ The word $\omega $ is a
rop-word if and only if the two conditions are satisfied:

(i) The length of $\omega $ is odd, say $2\ell +1$.

(ii) The height of the prefixes of length $\ell $ of the rotations of $%
\omega $\ are equal to $h-2$ or $h-1$.
\end{lemma}

From this lemma, Andr\'{e} Bouchet shows the following theorem. Let $\omega $
be a word of length $n$ and $d$ be an integer such that $0 <d \leq n/2$. A $%
d $-pairing of $\omega $\ is a partition of the subset of indices $\{i :0
\leq i <n ,\omega _{i} =3\}$ into pairs of indices $\{j ,j +d\}$. Arithmetic
operations on indices are to be understood $\mathop{\rm mod} ~n$.

\begin{theorem}
Let $\omega $ be a word of even height.\ $\omega $ is a rop-word if and only
if the two conditions are satisfied:

(i) The length of $\omega $ is odd, say $2 \ell +1$.

(ii) $\omega $ admits a $\ell $-pairing.
\end{theorem}

\begin{figure}[h]
\centering
\includegraphics[width=7cm]{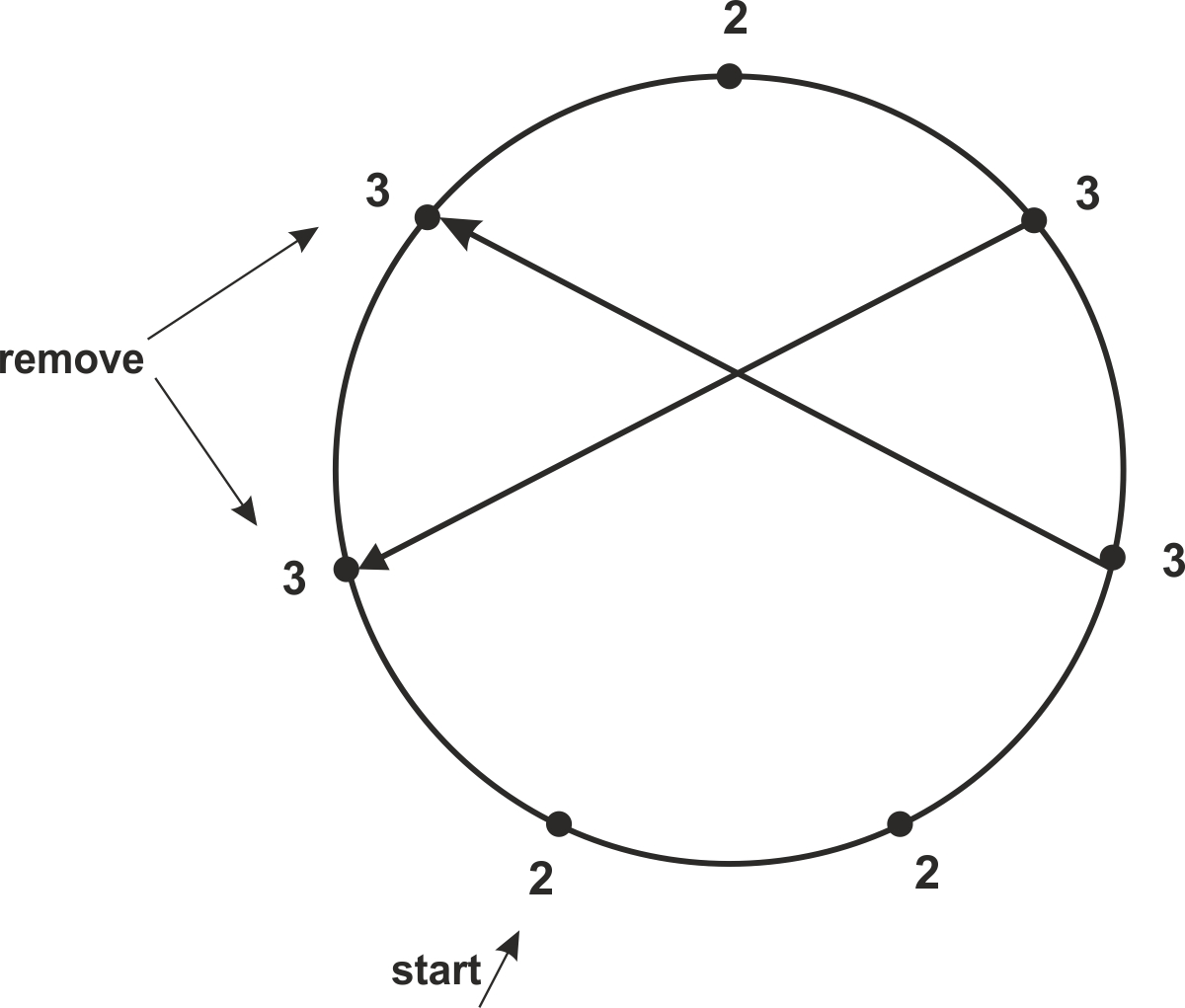}
\caption{Cyclic representation of rop-words}
\label{Cycle}
\end{figure}

Let $n_{2}$ and $n_{3}$\ be the number of symbols 2 and 3 in $\omega $ and $%
n =n_{2} +n_{3}$ be the length of $\omega \text{.}$ For a given $n_{2}$, we
will use the $d$-pairing to show that there is a one to one correspondance
between aperiodic necklaces of length $n$ with $n_{2}$ black beads
(represented by letter 2) and $(n -n_{2})$ white beads (represented by
letter 3) and Lyndon rop-words of length $n^{ \prime } =2 n -n_{2}$ with $%
n_{2}^{ \prime } =n_{2}$ letters 2 and $n_{3}^{ \prime } =2 (n -n_{2})$
letters 3. And also a one-to-one correspondance between necklaces
(eventually periodic) of length $n$ with $n_{2}$ black beads (represented by
letter 2) and $(n -n_{2})$ white beads (represented by letter 3) and
rop-words.\ The correspondance is obtained by adding or removing the letters
3 coming from the pairing.\ Let us examine an example. (See fig. \ref{Cycle})

Fix $n_{2}$, for instance $n_{2}=3$, and let $n$ be $n=5$. The word 2233233
is a (Lyndon) rop-word with odd length 7 ($n_{2}^{\prime }=3$, $%
n_{3}^{\prime }=4$) since it has a 3-pairing.\ Put the word on a circle,
starting from the bottom and turn counterclockwise as shown on the figure %
\ref{Cycle}. Now discard the second 3 of each pairing $(3,3)$ turning
counterclockwise. Reading the remaining word clockwise starting from the
bottom gives the word 22332, one of the two necklaces of length 5 with 3
letters 2. Conversely, starting from the word 22233, it is easy to add a
3-pairing by doubling each letter 3, with respect to the counterclockwise
tour.

Let $\omega =\omega _{0}\omega _{1}\ldots \omega _{n-1}$ be a word of $%
\{2,3\}^{\ast }$ and $p$ coprime with $n$. Denote by $\omega
^{(p)}=x_{0}x_{1}\ldots x_{n-1}$ the word obtained by reading all letters of 
$\omega $ by step $p$, starting from $\omega _{0}\text{.}$ Namely, each
letter of $\omega ^{(p)}$\ is $x_{j}=\omega _{k}$ with $k=jp~\mathop{\rm mod}
$ $p$, $0\leq j<p\text{.}$ For instance, the word $\omega =2233233$ depicted
on fig. \ref{Cycle} with $n=7$ and $\ell =3$ becomes $\omega ^{(2)}=233233.$
A.\ Bouchet shows

\begin{theorem}
Let $\omega $ be a word of even height.\ $\omega $ is a rop-word if and only
if the two conditions are satisfied:

(i) The length of $\omega $ is odd, say $2 \ell +1$.

(ii) $\omega ^{(\ell )}$ admits a $1$-pairing.
\end{theorem}

In other words, we can always transform a (resp. Lyndon) rop-word $\omega $
of length $2 \ell +1$\ by a one-to-one map $\phi $\ such that the letters 3
in $\omega ^{(\ell )}$ are always coupled by subwords 33.\ The bijection $%
\psi $ sending $2 \rightarrow 0$\ and\ $33 \rightarrow 1$ maps $\omega
^{(\ell )}$ to a word $\omega ^{ \prime } \in \{0 ,1\}^{ \ast }$
corresponding to a (resp. aperiodic) necklace. 
\begin{equation*}
\omega \overset{\phi }{ \longrightarrow }\omega ^{(\ell )}\overset{\psi }{
\longrightarrow }\omega ^{ \prime }
\end{equation*}
The table \ref{Iso} shows the first Lyndon rop-words for $n_{2} =3$ and the
corresponding aperiodic necklaces

\medskip 
\begin{table}[h]
\centering
\begin{tabular}{llll}
\hline
Aperiodic Necklaces & $n$ & Lyndon Rop-words & $n^{\prime }$ \\ \hline
\multicolumn{1}{c}{0001} & \multicolumn{1}{c}{4} & \multicolumn{1}{c}{22323}
& 5 \\ 
\multicolumn{1}{c}{00011} & \multicolumn{1}{c}{5} & \multicolumn{1}{c}{
2233233} & 7 \\ 
\multicolumn{1}{c}{00101} & \multicolumn{1}{c}{5} & \multicolumn{1}{c}{
2323233} & 7 \\ 
\multicolumn{1}{c}{000111} & \multicolumn{1}{c}{6} & \multicolumn{1}{c}{
223332333} & 9 \\ 
\multicolumn{1}{c}{001011} & \multicolumn{1}{c}{6} & \multicolumn{1}{c}{
232332333} & 9 \\ 
\multicolumn{1}{c}{001101} & \multicolumn{1}{c}{6} & \multicolumn{1}{c}{
232333233} & 9 \\ \hline
\end{tabular}
\medskip
\caption{Correspondance for $n_{2}=3$}
\label{Iso}
\end{table}

Conversely, starting from the representing Lyndon word of a aperiodic
necklace $\omega ^{ \prime }\text{,}$\ we construct the word $\omega ^{(\ell
)}$\ by the bijection $\psi ^{ -1}$ sending $0 \rightarrow 2$\ and\ $1
\rightarrow 33\text{,}$ and the word $\omega $ by applying $\phi ^{ -1}\text{%
.}$ By construction, the height $h (\omega ^{(\ell )})$\ is even and also $h
(\omega ) .\;$Moreover, $\omega $\ has a $\ell $-pairing and then is a
rop-word.

\section{Enumeration of the rop-words}

The number of necklaces (see \cite{Aig2007} for details, and also \cite%
{Cas2013} for applications to music theory) with $n_{2}$ black beads and $%
n_{3}/2$ white beads derives from the generating function of the action of
the cyclic group 
\begin{equation}
Z(C_{n_{2}},x)=\frac{1}{n_{2}}\sum_{d|n_{2}}\varphi (d)x_{d}^{n_{2}/d}
\end{equation}%
where the sum is over all divisors $d$ of $n_{2}$ and $\varphi $ is the
Euler totient function, according to the substitution of $x_{j}$ by $\frac{1%
}{1-x^{j}}\text{.}$ The developpement gives the coefficients of $x_{j}$
which are precisely the number of necklaces with $n_{2}$ black beads and $j$
white beads. For example, if $n_{2}=p$ is prime, the developpement leads to
the following equations: 
\begin{eqnarray}
Z(C_{p},x) &=&\frac{1}{p}\left( \varphi (1)x_{1}^{p}+\varphi (p)x_{p}\right)
\\
&=&\frac{1}{p}\frac{1}{(1-x)^{p}}+\frac{p-1}{p}\frac{1}{1-x^{p}}  \notag \\
&=&\frac{1}{p}\left( 1+\overset{\infty }{\sum_{n=1}}\frac{p(p+1)\ldots
(p+n-1)}{n!}x^{n}\right) +\frac{p-1}{p}\left( \overset{\infty }{\sum_{n=0}}%
x^{np}\right)  \notag \\
&=&1+\frac{1}{p}\overset{\infty }{\sum_{n=1}}\binom{p+n-1}{n}x^{n}+\frac{p-1%
}{p}\left( x^{p}+x^{2p}+x^{3p}+\ldots \right)  \notag \\
&=&1+\overset{\infty }{\sum_{n=1}}a_{n}x^{n}  \notag
\end{eqnarray}%
with 
\begin{equation}
a_{n}=\left\{ 
\begin{array}{lll}
\binom{p+n-1}{n} & if & n%
\not%
\equiv 0~\mathop{\rm mod}~p \\ 
\binom{p+n-1}{n}+p-1 & if & n\equiv 0~\mathop{\rm mod}~p%
\end{array}%
\right.
\end{equation}%
The table \ref{Rn} with $n_{2}$ on the horizontal axis and $n_{3}$\ on the
vertical axis shows the number of rop-words for $n_{2}$ and $n_{3}$ fixed.
The number of rop-words of length $n$ is given by summing along the diagonal 
$n_{2}+n_{3}=n$. Each column of the table is obtained from the developpement
of the generating function $Z(C_{n_{2}},x)$.\ For $n_{2}$ prime, the
coefficients agree with the formula of $a_{n}$ given above.~

\begin{table}[h]
\centering%
\begin{tabular}{cccccccccc}
\hline
& 1 & 3 & 5 & 7 & 9 & 11 & 13 & 15 & 17 \\ \hline
2 & 1 & 1 & 1 & 1 & 1 & 1 & 1 & 1 & 1 \\ 
4 & 1 & 2 & 3 & 4 & 5 & 6 & 7 & 8 & 9 \\ 
6 & 1 & \emph{4} & 7 & 12 & \emph{19} & 26 & 35 & \emph{46} & 57 \\ 
8 & 1 & 5 & 14 & 30 & 55 & 91 & 140 & 204 & 285 \\ 
10 & 1 & 7 & \emph{26} & 66 & 143 & 273 & 476 & \emph{776} & 1197 \\ 
12 & 1 & \emph{10} & 42 & 132 & \emph{335} & 728 & 1428 & \emph{2586} & 4389
\\ \hline
\end{tabular}
\medskip
\caption{Number of rop-words $(n_{2},n_{3})$}
\label{Rn}
\end{table}
In each column, we recover the number of binary necklaces with length $%
n_{2}+q$ and density $q=n_{3}/2$\ given by the rhs of the next formula. From
the bijection of the previous section, it follows that~the number $%
R(n_{2},n_{3})$ of rop-words with $n_{2}$ symbols 2 and $n_{3}$ symbols 3 is
the number of binary necklaces of length $n_{2}+n_{3}/2$ and density $%
n_{3}/2 $,

\begin{equation}
R(n_{2},2q)=\frac{1}{n_{2}+q}\sum_{d|\gcd (n_{2}+q,q)}\varphi (d)\left( 
\begin{array}{c}
(n_{2}+q)/d \\ 
q/d%
\end{array}%
\right) ,\text{\quad \quad }q=1,2,3,\ldots
\end{equation}

By computing Lyndon words on~ alphabet~$\{2,3\}$ and deleting those which
are not rop-words, we get the table \ref{Ln}\ of the number of Lyndon
rop-words for $n_{2}$ and $n_{3}$ fixed, with $n_{2}$ on the horizontal axis
and $n_{3}$\ on the vertical axis. The total number of Lyndon rop-words of
length $n$ is obtained by summing along $n_{2}+n_{3}=n$. The differences
between the tables~\ref{Rn} and~\ref{Ln} are in italics.

\begin{table}[h]
\centering
\begin{tabular}{cccccccccc}
\hline
& 1 & 3 & 5 & 7 & 9 & 11 & 13 & 15 & 17 \\ \hline
2 & 1 & 1 & 1 & 1 & 1 & 1 & 1 & 1 & 1 \\ 
4 & 1 & 2 & 3 & 4 & 5 & 6 & 7 & 8 & 9 \\ 
6 & 1 & \emph{3} & 7 & 12 & 18 & 26 & 35 & \emph{45} & 57 \\ 
8 & 1 & 5 & 14 & 30 & 55 & 91 & 140 & 204 & 285 \\ 
10 & 1 & 7 & \emph{25} & 66 & 143 & 273 & 476 & \emph{775} & 1197 \\ 
12 & 1 & \emph{9} & 42 & 132 & 333 & 728 & 1428 & \emph{2584} & 4389 \\ 
\hline
\end{tabular}
\medskip
\caption{Number of Lyndon rop-words $(n_{2},n_{3})$}
\label{Ln}
\end{table}
In each column of table \ref{Ln}, we recover the number of fixed density
Lyndon words given by the following formula, with $n_{3}=2q$. It follows
from the previous section, that the number $L(n_{2},n_{3})$ of Lyndon
rop-words with $n_{2}$ symbols 2 and $n_{3}$ symbols 3 is 
\begin{equation}
L(n_{2},2q)=\frac{1}{n_{2}+q}\sum_{d|\gcd (n_{2}+q,q)}\mu (d)\left( 
\begin{array}{c}
(n_{2}+q)/d \\ 
q/d%
\end{array}%
\right) ,\text{\quad \quad }q=1,2,3,\ldots
\end{equation}%
where $\mu $ is the Mobius function.

By summing these formulas along a diagonal $n=n_{2}+n_{3}$, we get the
number $\mathcal{L}_{n}$\ of Lyndon rop-words of length $n$\ and the number $%
\mathcal{R}_{n}$\ of rop-words of the length $n$:

\begin{equation}
\mathcal{L}_{n}=\sum_{n_{2}+n_{3}=n}L(n_{2},n_{3})=\overset{(n-3)/2}{%
\sum_{p=0}}~L(2p+1,n-2q-1)
\end{equation}%
and 
\begin{equation}
\mathcal{R}_{n}=\sum_{n_{2}+n_{3}=n}R(n_{2},n_{3})=1+\overset{(n-3)/2}{%
\sum_{p=0}}~R(2p+1,n-2p-1)
\end{equation}%
These numbers are tabulated as follows: 
\begin{table}[h]
\centering%
\begin{tabular}{cccccccccccccc}
\hline
$n$ & 3 & 5 & 7 & 9 & 11 & 13 & 15 & 17 & 19 & 21 & 23 & 25 & 27 \\ 
$\mathcal{R}_{n}$ & 2 & 3 & 5 & 10 & 19 & 41 & 94 & 211 & 493 & 1170 & 2787
& 6713 & 16274 \\ 
$\mathcal{L}_{n}$ & 1 & 2 & 4 & 8 & 18 & 40 & 90 & 210 & 492 & 1164 & 2786 & 
6710 & 16264 \\ \hline
\end{tabular}
\medskip
\caption{Numbers of Lyndon rop-words and rop-words of length $n$}
\label{Tot}
\end{table}
If $n$\ is prime, the difference between the cardinal of the two sets is 1,
since the word $2^{n}$ (where the letter 2 is repeated $n$ times) is a
rop-word but not a Lyndon word. If $n$ is a product or a power of primes,
some periodic words appear that are rop-words but not Lyndon words.\ This
explains the differences between\ the set of rop-words and the set of \emph{%
Lyndon} rop-words. For instance, if $n=9$, $(233)^{3}$ is a rop-word but not
a Lyndon word. The same is true for the words $(22323)^{3}$, $(233)^{5}$ and 
$(23333)^{3}$ of length 15. For $n_{2}=9$ and $n_{3}=12$, there are 333
Lyndon rop-words and 335 rop-words.\ The two non Lyndon rop-words are: $%
(2233233)^{3}$ and $(2323233)^{3}$.

\section{Generalization: from rop to rap words}

\subsection{First generalization.\ }

Let $s$ be an integer $\geq 2.$ The rhythmic oddity property could be
generalized in two ways.\ The first way is as follows.

\begin{definition}
A word $\omega \in \{2 ,3\}^{ \ast }$ is a $s$-rop word if

(i) $h (\omega ) \equiv 0~\mathop{\rm mod} ~s$

(ii) No cyclic shift of $\omega $ can be factorized into $s$ words $u_{1}
,u_{2} ,\ldots ,u_{s}$ such that 
\begin{equation*}
h (u_{1}) =h (u_{2}) =\ldots =h (u_{s})
\end{equation*}
\end{definition}

A Lyndon $s$-rop word is both a $s$-rop word and a Lyndon word.\ For
instance, if $s=3$, the first 3-rop words are: 2223, (3333), 22233, 22323,
(33333), 222333, 2222223, 2232333, 2233233, 2233323, 2323233, (3333333).\
Non Lyndon words are given in parenthesis. A computation of the number of
the first 3-rop words of length $n$ is given in table \ref{T5}.

\begin{table}[h]
\centering
\begin{tabular}{ccccccccccccccc}
\hline
$n$ & 4 & 5 & 6 & 7 & 8 & 9 & 10 & 11 & 12 & 13 & 14 & 15 & 16 & 17 \\ 
$\mathcal{R}_{n}^{(3)}$ & 2 & 3 & 1 & 6 & 11 & 6 & 25 & 46 & 41 & 117 & 232
& 278 & 631 & 1237 \\ 
$\mathcal{L}_{n}^{(3)}$ & 1 & 2 & 1 & 5 & 9 & 6 & 22 & 45 & 40 & 116 & 226 & 
278 & 620 & 1236 \\ \hline
\end{tabular}%
\medskip
\caption{Number of 3-rop words}
\label{T5}
\end{table}
For $s=4$, the first 4-rop words are: 22233, 22323$^{\ast }$, (222222),
223333, 232333, (233233), 2222233, 2222323, 2223223$^{\ast }$, 2333333$%
^{\ast }$. Non Lyndon words are in parenthesis.\ The star indicates 2-rop
words. A computation of the number of the first 4-rop words of length $n$ is
given in table \ref{T6}

\begin{table}[h]
\centering
\begin{tabular}{cccccccccccccc}
\hline
$n$ & 5 & 6 & 7 & 8 & 9 & 10 & 11 & 12 & 13 & 14 & 15 & 16 & 17 \\ 
$\mathcal{R}_{n}^{(4)}$ & 2 & 4 & 4 & 5 & 13 & 27 & 47 & 50 & 131 & 284 & 479
& 685 & 1450 \\ 
$\mathcal{L}_{n}^{(4)}$ & 2 & 2 & 4 & 5 & 12 & 24 & 47 & 50 & 131 & 279 & 473
& 683 & 1440 \\ \hline
\end{tabular}%
\medskip
\caption{Number of 4-rop words}
\label{T6}
\end{table}

\subsection{Second generalization.}

The second way of generalization is to change the alphabet and to consider
words over $\mathcal{A} =\{1 ,2 ,\ldots ,s\}$.

\begin{definition}
A word $\omega \in \{1 ,2 ,\ldots ,s\}^{ \ast }$ has the rhythmic arity
property (rap) of order s if

(i) $h (\omega ) \equiv 0~\mathop{\rm mod} ~s$

(ii) No cyclic shift of $\omega $ can be factorized into $s$ non-empty words 
$u_{1} ,u_{2} ,\ldots ,u_{s}$ such that 
\begin{equation*}
h (u_{1}) =h (u_{2}) =\ldots =h (u_{s})
\end{equation*}
\end{definition}

For short, we call \emph{s-rap-word} a word with the rhythmic arity property
of order $s$.\ For example, on the alphabet \{1,2,3\}, the words 111 and
333\ are not 3-rap-words, but 123 and 132 are. The word 11133 is not a
3-rap-word since the subwords $u_{1}=111$, $u_{2}=3$ and $u_{3}=3$ have the
same height, but the word 11313 is a rap-word. A\ \emph{Lyndon s-rap-word}
is both a $s$-rap-word and a Lyndon word.\ For instance, 11133, 11313, 11322
are Lyndon 3-rap-words.

\section*{Acknowledgements}

\thispagestyle{empty} I would like to thank Harald Fripertinger for valuable
comments and remarks.

\end{document}